\documentclass[aps,amsmath,amssymb,floatfix,superscriptaddress,nofootinbib]{revtex4}

\usepackage{graphicx}
\usepackage{dcolumn}
\usepackage{bm}
\usepackage{latexsym}
\usepackage{latexsym,graphics}

\newcommand{\be}{\begin{equation}}
\newcommand{\ee}{\end{equation}}
\newtheorem{lemma}{Lemma}

\begin{document}

\title{k-workers in a Circular Warehouse -- A random walk on a circle, without passing.}

\author{Joseph D. Skufca}
\email{skufca@usna.edu}
\affiliation{University of Maryland, College Park}
\altaffiliation{United States Naval Academy, Annapolis, Maryland}

\date{\today}

\begin{abstract} We consider the problem of stochastic flow of multiple particles traveling on a closed loop, with a constraint that particles move without passing.  We use a Markov chain description that reduces the problem to a generalized random walk on a hyperplane (with boundaries).  By expressing positions via a moving reference frame, the geometry of the no-passing criteria is greatly simplified, with the resultant condition expressible as the coordinate system planes which bound the first orthant. To determine state transition probabilities, we decompose transitions into independent events and construct a digraph representation in which calculating transition probability is reduced to a shortest path determination on the digraph.  The resultant decomposition digraph is self-converse, and we exploit that property to establish the necessary symmetries to find the stationary density for the process.

\end{abstract}

\maketitle

\section{Introduction}

In dense traffic flow, the speed of individual vehicles is constrained because the cars are unable to pass one another.  In protein controlled diffusion through a membrane, the narrow passage forces the transported molecules to move single-file.  A similar phenomena develops in the motion of people along multiple serving stations, such as occurs in a buffet line.  Motivated by the natural problems that develop when physical objects try to move independently along the same path, we consider a stochastic flow of particles on a closed loop, with the constraint that one particle may not pass another.  We limit this paper to the specific problem where motion occurs over discrete units of time, and the particles are limited to move in some integer number of steps.  We expect that our approach may have some applicability to models of microscopic transport, cellular automata models of traffic flow, and some queuing theory problems, as well as providing some general insight on diffusive flows in constrained channels. 

{\bf Problem Statement.} Consider a warehouse with parts stored in $n$ locations (parts bins) arranged in a closed path.  $k$ workers move around the warehouse in a continuous loop, filling an order on each trip through the warehouse.  In real life, a worker would collect from those bins which contain the parts required to fill an order.  We make a simplifying assumption and model this selection process as a Bernoulli trial, with probability $s$ that the worker will select a part from any given bin.  All the workers move in the same direction as they fill their orders, but without passing each other.  (Perhaps they are all driving forklifts, and there is no room to pass; or maybe the workers are actually robotic carts moving on a common track.)  To define ``time," we assume that during one time step, each worker moves around the closed loop until one of two conditions is met: \begin{enumerate}
\item The worker has a Bernoulli success at some bin and stops to collect a part, or 
\item The worker is ``blocked'' from proceeding around the loop by some other worker who has stopped to gather a part.  Since passing is not allowed, the first worker does not have a Bernoulli success during that time step and collects no part. \end{enumerate} Blockages, therefore, represent reduced output from a worker. Our goal is to evaluate the efficiency of our warehouse by {\it determining the expected value of the fraction of the time steps that an individual worker is blocked.}  

We apply the theory of Markov chains to this problem, which leads to a description of the process as a $k-1$ dimensional generalized random walk on a compact lattice.  We develop the state space description in a moving reference frame, where the constraint of single file motion can be simplified to a boundary condition expressible as a coordinate plane.  To determine transition probabilities, we develop an associated digraph that allows us to decompose transitions in a way that relates transition probabilities to shortest paths between vertices of that digraph.  We then exploit the self-converse structure of that digraph to establish a stationary density on the state space.  The paper is organized as follows: in \S  \ref{sstatespace}, we define and enumerate the states; in \S  \ref{stransitions}, we develop an expression for the transition probabilities; in \S \ref{sdistribution} we establish a stationary distribution for the general problem; and in \S \ref{sconclusion} we answer the specific efficiency question of the problem statement.

\section{Labeling and Counting States. \label{sstatespace} }

\begin{figure}[htbp]
\includegraphics{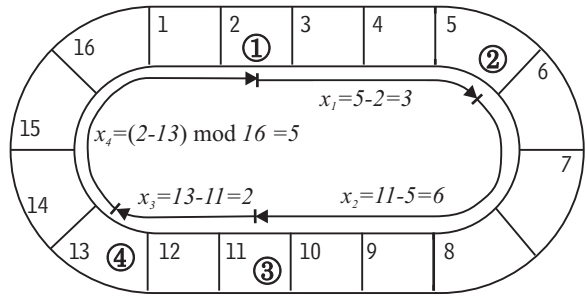} 
\caption{A warehouse of $n=16$ bins.  The positions of the $k=4$ workers are shown in circles. Note that $x_1+x_2+x_3+x_4=3+6+2+5=16=n.$  The bin locations of the four workers in the fixed frame, $(2,5,11,13)$ is a member of the equivalence class denoted by configuration $(3,6,2,5).$}
\label{fdistdef}
\end{figure}

Markov Chain theory requires that the probabilities of future states of the system depend only upon the current state of the system, so the initial step in applying these techniques is to identify a state space description that satisfies this requirement.  Since all the bins have the same probability for selection, a blockage event does not depend upon the specific bin location of the worker, but rather it is the distance from one worker to the next that is of critical importance.  Therefore, we focus on these distances to provide a state space description. Let
\[
d=w-z \pmod{n}
\]
be the distance from a worker at bin $z$ forward to a worker at bin $w$ (in the direction of motion of the workers).  We specify that two different workers may not occupy the same bin, so $d \geq 1.$  We let $X=(x_1,...,x_k),\  X \in \mathsf{R}^k,$ with $x_i$ the distance from worker $i$ forward to worker $i+1\pmod{n},$ and we say that vector $X$ defines a {\it physical configuration} of the system (or simply {\it configuration}).  This vector description can be viewed as using a moving reference frame (moving with the workers), which replaces the fixed frame of the bin numbers.  Each configuration is actually an equivalence class of bin locations for the workers.  Fig \ref{fdistdef} provides a visual example of this method of configuration identification.  

\begin{figure}[htbp]
\includegraphics{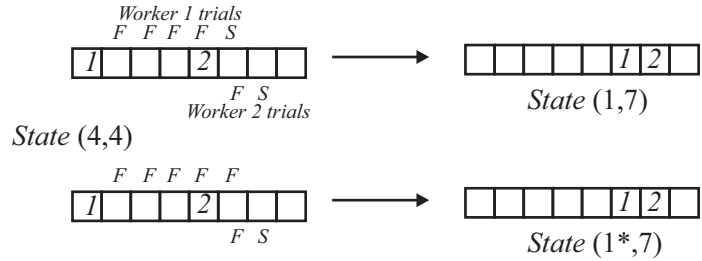} 
\caption{With $k=2$ and $n=8,$ we see two possible transitions from the state $(4,4).$ The Bernoulli trial results are indicated with $F$ for failure and $S$ for success.  Both transitions end in the same physical configuration.  However, in the second case, worker 1 does not have a Bernoulli success before he is forced to stop because he cannot occupy the same bin as worker 2.  In this second case, worker 1 would not collect a part.}
\label{fblocks}
\end{figure}

Since the workers are restricted to a warehouse of $n$ bins, we have $\sum_{i=1}^k x_i=n.$ Counting configurations is simply the problem: {\it how many ways can you arrange $k$ positive integers that sum to $n$?} This well known problem, generally covered in a first course in combinatorics, yields the formula
\be
\label{ecountconfig}
C(k,n)=  {n-1 \choose k-1}
\ee
possible {\it configurations.}  To extend the set of configurations to describe all the states of the system, we must consider blockages.  In a {\it blocked}  situation, the blocked worker is physically a  distance $d=1$ from the worker ahead.  The physical configuration is the same in both cases.  The difference, however, is that at the bin where he is stopped, the blocked worker had a Bernoulli failure instead of a Bernoulli success (and therefore does not collect a part on that time step).  To distinguish this condition, we use the symbol $1^*$ to represent that blocked condition. Fig \ref{fblocks} illustrates this difference. 

By defining $|1^*|:=1,$  we can formally represent states as vectors:
\be
\label{estates}
X=(x_1, \ldots ,x_k) \qquad x_i \in \{1^*,1,2,\dots,n-k+1\},\ \sum_{i=1}^k |x_i|=n.
\ee
Geometrically, the blockage condition forms a boundary on the configuration space, where these boundaries can be treated as hyperplanes that are parallel to each of the coordinate planes. \footnote{Alternatively, we could have used $0$ instead of $1*$ to represent a blockage, which makes it easier to interpret the blockage states as the coordinate hyperplanes $x_j=0$ forming the boundary of the physical configurations which lie strictly in the first orthant.  We find the $|1*|=1$ symbology more useful in providing an algebraic and physical description of the state.}  The simplicity of this boundary expression results from our having shifted to the moving coordinate system.  In a fixed reference (with states described by bin locations), the geometry of the blockage state becomes significantly more complex as $k$ increases \cite{BenAvraham}.

We enumerate the blocked states as follows: Suppose we want to count the number of states with $b$ blockages. An arbitrary state can be considered an ordered list of $k$ elements, and we seek to count the number of valid ways to fill those $k$ elements.  Forming such a sequence can be broken down into two steps: first, place a $1*$ for each worker that is blocked; and then fill the other elements with positive integers to satisfy the summation constraint of \eqref{estates}.  There are ${k \choose b}$ ways to insert the $1^*$ s into the stream.  The other $k-b$ elements must by filled with positive integers that sum to $n-b.$ We use \eqref{ecountconfig} with new arguments $(k-b)$ and $(n-b)$ to count the possibilities for the other integers, which is ${n-b-1 \choose k-b-1}$.  Applying the multiplicative rule of counting, we find $N,$ the number of states with $b$ blockages, is given by
\be
\label{nb}
N(b,k,n)={k \choose b}{n-b-1 \choose k-b-1}.
\ee

\vspace{.2in}
\noindent {\bf Example.} Let $k=5$ and $n=16.$  An arbitrary state is $(\square, \square, \square, \square, \square).$ Suppose $b=2.$ We must place two $1*$s into five possible slots; there are ${5 \choose 2}$ ways to do so, and   $(\square, \square, 1*, \square, 1*)$ would be a specific choice.  We now must fill the other three positions with positive integers which sum to $14,$ and there are ${13 \choose 2}$ ways to do so.  Then $N(2,5,16)={5 \choose 2}{13 \choose 2}=780.$
\vspace{.05in}

To count the entire state space, we simply sum over all possible values for $b.$  (We note that the maximum number of blockages is $b=k-1,$ since one worker must be stopped in order to block the other workers.) The total number of states would be given by
\be
\label{numstates}
M(k,n)= \sum_{b=0}^{k-1} N(b,k,n)= \sum_{b=0}^{k-1}{k \choose b}{n-b-1 \choose k-b-1}.
\ee

\section{Transition Probabilities \label{stransitions}}

Given the state description of \S  \ref{sstatespace}, we number those states from $1$ to $M.$ To establish the Markov property, we must identify the transition probabilities
\be
p_{ij}={\mathbb P} \left[ \mbox{System in state } X_j \mbox{ at time } \tau+1| \mbox{ System in state } X_i \mbox{ at time } \tau \right].
\ee
(For our general problem, there is no obvious scheme for a preferred numbering of the states.  Therefore, we develop transitions probabilities as $p_{XY},$ the probability of transitioning from state $X$ to state $Y,$ and identify $p_{ij} \equiv p_{X_iX_j}.$)  First, we will consider only the transitions between unblocked states, and then we extend that formulation to describe transitions from and to the blocked states.

\subsection{Transition probabilities for unblocked states.}

Recall that $s$ is the probability of Bernoulli success, and let $q=1-s.$ We decompose a transition from unblocked state $X$ to unblocked state $Y$ into a sequence of independent events:

\[
\begin{array}{lccc}
\mbox{Event:} & \underbrace{\parbox[b]{1.3in}{Each worker has $Q$ bernoulli failures {\it before} one of the workers has a success.}} & \underbrace{\parbox[b]{1in} { All workers have one success.}} & \underbrace{\parbox[b]{2.3in}{ After one worker has stopped, the other workers have additional failures to move the system from configuration $X$ to configuration $Y.$}} \\
\mbox{Probability:} & q^{kQ} & s^k & f(X,Y) 
\end{array}
\]
We compute the transition probability by taking the product of the independent events and summing over all $Q$ to give
\be
p_{XY}=\sum_{Q=0}^\infty q^{kQ}s^k f(X,Y)=s^k \frac{1}{1-q^k} f(X,Y),
\ee
where $f(X,Y)$ remains to be specified.  We call $f(X,Y)$ the {\it arrangement probability,} and note it is a conditional probability (conditioned upon one of the workers having stopped).  Our visualization is that all workers simultaneously move to an adjacent bin and conduct a trial. If all workers have failures, they would all index by one bin and repeat.  Until the first success, the physical configuration does not change. However, when one of the workers has a success, that worker stops, and the other workers keep moving.    The change in configuration is the result of workers moving -- Bernoulli failures.  Therefore, $f(X,Y)=q^l,$ where $l$ is the total number of bins that the workers must pass to take the system from configuration $X$ to configuration $Y,$ {\it given} that one of the workers passes 0 bins.  We use a geometric and graph theoretic approach to resolve the problem of counting these Bernoulli failures.

\begin{figure}[htb]
\vspace{2in}
\special{eps: 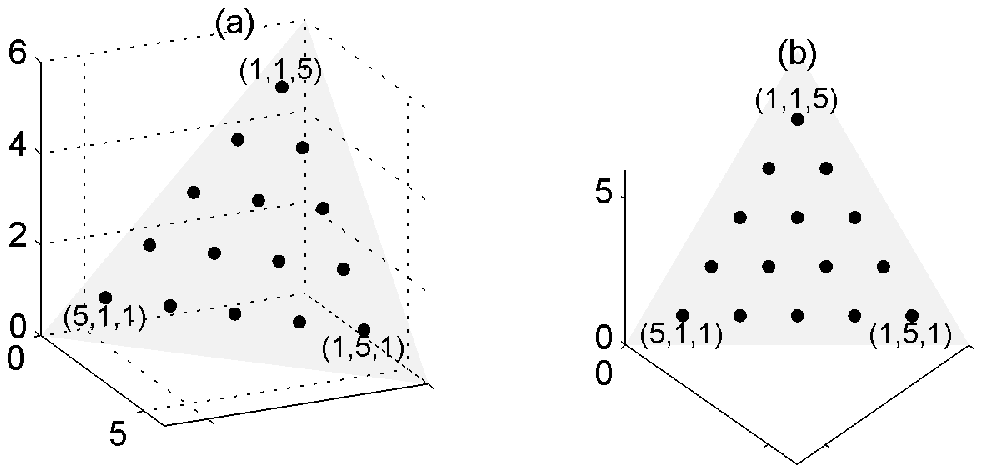 x=5in y=2in}
\caption{Visualization of states for $k=3, n=7.$}
\label{fig3d}
\end{figure}

{\bf Background geometry to address arrangement probability.} Let $X$ be an unblocked state -- a $k$-tuple of positive integers which sum to $n.$  Then $X$ lies on the hyperplane in $\mathbb{R}^k$ given by $\sum_{i=1}^k x_i = n,$ with orthogonal vector $u=(1,\ldots,1) \in \mathbb{R}^k.$  Specifically, $X$ is a lattice point in the first orthant that lies on that hyperplane. Define vector $V^r=(v_1,\ldots,v_k),$ by
\be
v_i=\begin{cases}
n-k+1& i=r ,\\
1& i \neq r.
\end{cases}
\ee
Geometrically, each $V^r$ specifies a vertex of the polytope containing all the configurations. Fig \ref{fig3d} illustrates this geometry for the case $k=3,n=7$: Fig \ref{fig3d}a shows how the hyperplane and lattice points are situated in $\mathbb{R}^3,$ and Fig \ref{fig3d}b shows the view from a position orthogonal to the plane, illustrating the structure on the in $k-1=2$ dimensional subspace. Vertices $V^i$ are labeled on both figures.

We define $\Omega$ to be the convex hull of $\{V^r:r=1,\ldots,k\}.$ For $X,Y \in \Omega,$ we define $\Delta_{XY} = (d_1,\ldots,d_k):=Y-X.$ Since $\sum x_i=\sum y_i=n,$ a necessary condition for $\Delta$ is
\be
\label{dreq}
\sum_{i=1}^k d_i=0.
\ee
Geometrically, each $\Delta_{XY}$ is a lattice point on the hyperplane orthogonal to vector $u=(1,\ldots,1)$ that passes through the origin.  We call this hyperplane $\mathcal{P}.$  We seek to describe $f(X,Y)$ in terms of some function $g(\Delta_{XY}),$ with $g(\cdot)$ independent of $X$ and $Y.$ The justification for this approach is that if no blockages occur, the workers are not affected by each other. The rearrangement probability can be completely characterized by the {\it change} in spacing between workers, independent of the initial and final configurations.  If such function $g$ exists, then the process can be treated as a generalized random walk on ${\mathbb R}^{k-1}$ (the hyperplane dimension). The physical configurations simply characterize which $\Delta$ vectors are feasible, because we require $X+\Delta = Y \in \Omega.$ Since $X$, $Y$, and $\Omega$  are contained within a bounded set, our description becomes a random walk with boundary.

\begin{figure}[htb]
\vspace{1in}
\special{eps: 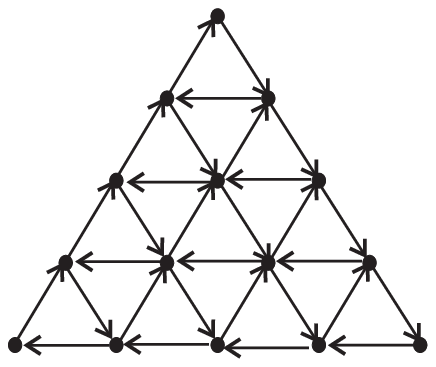 x=1in y=1in}
\caption{A digraph representation of rearrangements for $k=3$ and $n=7$ using $\delta$ vectors.  Each edge represents the effect of a single Bernoulli failure by one worker.}
\label{digraph}
\end{figure}

{\bf Canonical representation of $\Delta.$} Suppose $X$ describes a configuration. If worker $i+1$ has a Bernoulli failure and moves one bin, then the distance from the worker behind has increased by one, while the distance to the worker ahead has decreased by 1.  The new configuration $\tilde{X}$ would be given by
\[
\left\{\begin{array}{ll}
\tilde{x}_i=x_i+1 & \mbox{} \\
\tilde{x}_{i+1}=x_{i+1}-1 & \mbox{}\\
\tilde{x}_j=x_j & j \notin \{i,i+1\} \end{array} \right..
\]
Motivated by this effect of a single Bernoulli failure, we define vectors $\delta_i =(\delta_{i1},\ldots,\delta_{ik}),$ where
\[
\delta_{ij}= \left\{ \begin{array}{ll}
1 & j=i \\
-1 & j=i+1 \pmod{n} \\
0 & \mbox{else} \end{array} \right..
\]
[Example: $\delta_3=(0,0,1,-1,0,\ldots,0).$] We note that $\left\{ \delta_i \right\}_{i=1}^k \subset \mathcal{P},$ and $\mbox{span}\left( \left\{ \delta_i \right\}_{i=1}^k \right)=\mathcal{P}.$   These spanning vectors have a natural representation as edges in a digraph which connect the configurations that result from a single Bernoulli failure. We call this representation the {\it rearrangement digraph.}  Fig (\ref{digraph}) shows the rearrangement digraph for $k=3,n=7.$   

Any point $\Delta \in \mathcal{P}$ is a linear combination of the $\delta_i$s, so there exists vector $\alpha$ with integer components such that
\be
\label{deltaalpha}
\Delta=\sum_{i=1}^k \alpha_i \delta_i.
\ee
Since $\mathcal{P}$ is $k-1$ dimensional, the representation is not unique.  However, we want that decomposition to satisfy some additional constraints:  As previously stated, rearrangement probability is conditioned upon one of the workers having no failures.  Therefore, $\alpha_j=0$ for some $j \in \{1,\ldots,k\}.$  There are at most $k$ representations that satisfy this requirement.  We also need the representation to be a path in the rearrangement digraph, which implies $\alpha_i \geq 0$ for all $i.$  The following paragraphs describe the algorithm that identifies this canonical representation.
 
{\bf Construction of canonical representation.} Define matrix $A_{k \times k} = \left[\ \delta_1,\ldots,\delta_k \right],$ where the $\delta_i$s form columns of $A.$ Solutions to (\ref{deltaalpha}) are equivalent to solutions of $A\alpha=\Delta.$ We note that $\mbox{rank }A=k-1,$ so $A$ is singular.  It is easy to verify that the eigenvector for eigenvalue $\lambda=0$ is $u=(1,\ldots,1).$ The one parameter family of solutions to (\ref{deltaalpha}) is given by
\be
\label{afam}
\alpha=\alpha^0 + tu,
\ee
where $\alpha^0$ is any solution to (\ref{deltaalpha}) and $t$ an integer.

The geometric interpretation of $\alpha$ is that each $\alpha$ describes a path associated with the displacement $\Delta,$ and each $\alpha_i$ describes how many steps we take in the $\delta_i$ direction when traversing that path.  Since we want the decomposition to represent a valid rearrangement of the workers, we seek vector $\beta$ such that 
\be
\label{breq}
\left\{\begin{array}{ll}
A\beta=\Delta & \mbox{}, \\
\beta_i \geq 0 & \forall i \in \{1,\ldots,k\}, \mbox{ and}\\
\beta_j = 0 & \mbox{for at least one } j \in \{ 1,\ldots ,k \}. \end{array} \right\}.
\ee
To construct this vector $\beta,$ we proceed as follows: Define $\gamma := (\gamma_1,\ldots,\gamma_k),$ where 
\be
\label{egamma}
\gamma_j=\sum_{i=1}^{j} d_i.
\ee From \eqref{dreq} , we see that $\gamma_k=0,$ simplifying the computation which verifies that 
\be
A\gamma=\left[\begin{array}{c} \gamma_1-\gamma_k \\ -\gamma_1+\gamma_2 \\ \vdots \\ -\gamma_{j-1} + \gamma{j} \\ \vdots \end{array} \right]=\left[ \begin{array}{c} d_1+0 \\ -d_1+(d_1+d_2) \\ \vdots \\ -\sum_{i=1}^{j-1}d_i+\sum_{i=1}^j d_i \\ \vdots \end{array}\right]=\left[ \begin{array}{c} d_1 \\ d_2 \\ \vdots \\ d_j \\ \vdots \end{array} \right]=\Delta.
\ee
In general, $\gamma$ will not satisfy the second and third elements of \eqref{breq}.  To find the desired solution, we define $\gamma_{\min}=\mathop {\min }\limits_i \gamma_i,$ and assign $\beta_i=\gamma_i-\gamma_{\min}.$ Clearly, all of the components of $\beta$ are non-negative, with $\beta_j=0$ if $\gamma$ achieves its minimum at $j.$ In vector form, $\beta=\gamma-\gamma_{min}u.$  Since $A\gamma=\Delta,$ we use (\ref{afam}) to see that $\beta$ must be the solution to (\ref{breq}).  

\vspace{.2in}
\noindent{\bf Example:} Let $k=5$ and $n=12.$  Suppose $X=(2,1,4,4,1)$ and $Y=(3,1,2,2,4).$ Then $\Delta=Y-X=(1,0,-2,-2,3).$  We first find a solution to $A\gamma=\Delta$ by using \eqref{egamma} to find $\gamma=(1,1,-1,-3,0).$  Since $\gamma_{min}=-3,$ we find $\beta=(4,4,2,0,3).$ 
\vspace{.05in}

$\beta_i$ gives the number of Bernoulli failures required of worker $i+1$ to achieve the rearrangement specified by $\Delta,$ and $\|\beta\|=\sum \beta_i$ gives the total number of failures to achieve rearrangement. \footnote{In this paper, the only vector norm used is the 1-norm.  For ease of notation, we drop the subscript notation for vector norms.} Since any other rearrangement $\alpha=\beta+tu$ that describes $\Delta$ would result in either negative components of $\alpha,$ (for $t<0$) or all components positive (for $t>0$), both of which would fail to satisfy \eqref{breq},   the representation is unique.  Additionally, any other path representation would require $t>0$ and $\| \alpha \| >\|\beta\|.$ Therefore, $\| \beta \|$ describes the shortest path length.  We define $\phi (\Delta) = \|\beta\|,$ where $\beta$ is the unique vector satisfying (\ref{breq}).  We now can state that the rearrangement probability is given by
\[
f(X,Y)=g(\Delta_{XY})=q^{\left[\phi(\Delta_{XY})\right]}.
\]  
To simplify notation, we let $l_{ij}=  \phi(\Delta_{X_iX_j}).$  Then
\be
\label{pxy}
p_{ij}=s^k\frac{1}{1-q^k}q^{l_{ij}}
\ee
gives the transition probability {\it from an unblocked state to an unblocked state.}

\subsection{Transition Probabilities for Blocked States.}

Transition probabilities {\it from} blocked states is trivial: since the system does not care about past Bernoulli trials, transition probabilities depend only upon the initial configuration, not whether that configuration was the result of a blockage. So if we define a subtraction operation of $1^*$ by $m-1^*:=m-1,$ and calculate $\Delta$ using that arithmetic, transitions probabilities {\it from} blocked states are the same as the transition probability from the unblocked state of the same physical configuration.
 
To determine transition probabilities {\it to} blocked states, we note the following: when worker $j$ is {\it blocked,} his final action is the result of a Bernoulli failure, with no Bernoulli success.  So if we compare the two transitions from some state $Y,$
\be
\label{compblock}
\begin{array}{ll}
T_1: & Y\rightarrow(x_1,\ldots,x_{j-1},1,x_{j+1},\ldots,x_k) \\
T_2: & Y\rightarrow(x_1,\ldots,x_{j-1},1^*,x_{j+1},\ldots,x_k) \end{array},
\ee
we can express one transition probability in terms of the other. Specifically,
\be
\label{t2t1r}
P[T_2]=P[T_1]\cdot \frac{q}{s},
\ee
indicating that regardless of the specific collection of motions of the workers, transition $T_2$ requires one fewer success and one additional failure when compared to the $T_1$ transition. We apply \eqref{t2t1r} recursively to (\ref{pxy}) (once for each blockage) to yield:
\be
\label{pij}
p_{ij}=r^{b_j}\frac{s^k}{1-q^k}q^{l_{ij}} ,
\ee
where $r=q/s$ and $b_j$ is the number of blockages in state $X_j.$  This formula reduces to \eqref{pxy} when $b_j=0$, and therefore, it applies to all transitions.

\section{Stationary Distribution \label{sdistribution}}

A direct approach for finding stationary distributions is to develop a {\it state transition matrix,} $P=\left[ p_{ij} \right].$  Then the unit left eigenvector of $P$ associated with eigenvalue $\lambda=1$ gives a stationary distribution of the states\cite{Grimmett}. However, finding such eigenvectors for a general problem is difficult unless $P$ is very well structured.  Although we have a formula for $p_{ij},$ we can find no ``good" ordering of the states, where we interpret ``good" as giving enough structure to the state transition matrix to directly determine the eigenvector for the general problem of $k$ workers and $n$ bins.  

Often, finding general solutions requires some insight (in other words, a good guess).  Some of the observations that lead to a ``good guess" for this problem are: (1) since the probabilities are the same for each worker, we expect significant symmetry in the problem, (2) the $\Delta$ vector description of transitions yields a random walk, and we see no reason for the random walk to prefer any particular unblocked state, and (3) blocked states are closely related to the associated physical configuration, acting like boundaries on the state space, and that relationship is governed by the ratio $r=q/s$ and the number of blockages.  So motivated, we postulate an initial distribution and show that it is stationary.

{\bf Postulated density and outline of proof of stationarity.} Let $\{X_i\}_{i=1}^M$ be the collection of all the states and $b_i$ the number of blockages in state $X_i.$ Define 
\be
\label{eomega}
\omega_i=r^{b_i}.
\ee  
We will treat these as weighting values for the density in each state, and let
\be
\label{nubal}
\nu=\frac{1}{\|\omega\|}\omega
\ee
be the initial distribution of states.  {\bf We want to show} that $\nu$ is stationary, or equivalently, that $\nu_j=\sum_{i=1}^M \nu_ip_{ij}$ for all $j.$  One way to describe this balance is to say that at each time step, the density that transitions {\it from} state $X_j$ to all states is exactly equal to the sum of all the density that transition {\it into} $X_j.$  To establish the required balance, we will prove the following logical argument:
\begin{enumerate}
\item Regardless of the density at time $\tau,$ at time $\tau + 1,$ the density of a state $X$ with $b$ blockages is $r^b$ times the weight of the unblocked state which has the same physical configuration as $X.$ 
\item Each unblocked state can be placed in a one-to-one correspondence with another unblocked state, (which we call the {\it co-state}) that has the property that when the density is $\nu,$ the density that transitions {\it into} the state is the same as transitions {\it from} the co-state.  Since $\nu$ assigns the same density to all unblocked states, the implication is that $\nu$ is stationary on the unblocked states. 
\item Since the densities of the blocked states are completely determined by item (1) (which is satisfied by $\nu$), stationarity of the unblocked states immediately implies stationarity of the blocked states.
\end{enumerate}

The first element of the proof is strait-forward.

\begin{lemma}\label{lemtoblocks} [Item (1)] Let $W^b$ be a state with $b$ blockages and let $W$ be the state with the same physical arrange as $W^b,$ but with no blockages. If the total density that is transitioned into state $W$ is $\alpha,$ then the total density that transitions into state $W^b$ will be $r^b \alpha .$ \end{lemma}

\noindent {\bf Proof of lemma.}
From repeated application of (\ref{t2t1r}), we know that for every state $X,$ 
\be
\label{probr}
p_{XW^b}=r^bp_{XW}.
\ee
Suppose the density of state $X_i$ at time $\tau$ is $\xi_i.$ Let 
\[
\alpha=\sum_{i} \xi_{i} p_{X_iW}
\]
be the total density transitioning into $W.$ Then the density transitioning into $W^b$ is 
\be
\sum_{i} \xi_{i} p_{X_iW^b}=\sum_{i} \xi_{i} r^b p_{X_iW} =r^b \alpha. \qquad \Box
\ee
\vspace{.2in}

Item (3) of the proof follows directly from items (1) and (2), so the crux of the argument is rigorously proving that the unblocked states are stationary.  Although an algebraic argument can be shown, it is essentially a tedious symbolic manipulation that fails to reveal the underlying structure.  Instead, we will develop an proof that exploits (and illuminates) the symmetry that is visible in the graph theoretic description (the rearrangement digraph) previously introduced.

{\bf Motivating, defining, and analyzing {\it co-states}.}\label{counblocked}  We want to show that $\nu_j(\tau+1)=\nu_j(\tau)$; that the density of state $j$ is the same at time $\tau+1$ as at time $\tau.$  Direct computation is intractable.  Instead, we show that for each $X_j,$ we can match in with another indexed state, $X_{j^*},$  such that $\nu_j(\tau+1)=\nu_{j^*}(\tau).$  By assumption on $\nu,$ we know that if both states are unblocked, then $\nu_j(\tau)=\nu_{j^*}(\tau),$ and we achieve our goal (for unblocked states).  An alternate description is that we seek to find a state $X_{j^*}$ so that the density that the transitions out of state $X_{j^*}$ is the same as that transitioning into state $X_j.$  Algebraically, we want
\be
\label{baseco}
\sum_{i=1}^M\nu_{j^*} p_{j^* i} = \sum_{i=1}^M\nu_i p_{i j}.
\ee

\begin{figure}[htbp]
\includegraphics{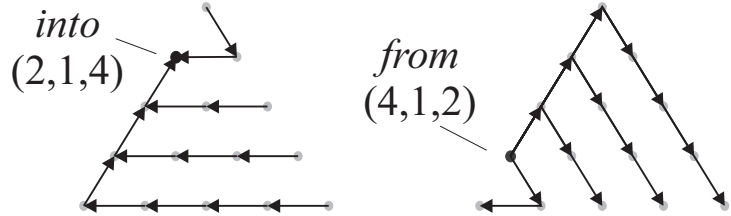} 
\caption{Using $k=3$ and $n=7,$ rearrangements into state $(2,1,4)$ and from its co-state $(4,1,2)$ are decomposed by $\delta$ vectors, to illustrate an equivalent set of pathlengths.}
\label{fcostate}
\end{figure}

We first consider a slightly different constraint than (\ref{baseco}).  Suppose we seek $j^*$ such that 
\be
\sum_{i=1}^{C}\nu_{j^*} p_{j^* i} = \sum_{i=1}^{C}\nu_i p_{i j},
\ee
where the summation is over only unblocked states.  Since $\nu_i=\nu_j$ for all unblocked states, we can use (\ref{pxy}) to simplify to 
\be
\label{cosum}
\sum_{i=1}^{C} q^{\left[\phi (\Delta_{X_{j^*} X_i})\right]} = \sum_{i=1}^{C} q^{\left[ \phi (\Delta_{X_i X_j})\right]}.
\ee
Since $q$ is arbitrary, \eqref{cosum} requires the equality of two polynomial expressions in $q,$ which implies that both the left and right sides must have the same collection of exponents.  In the digraph representation, each exponent ($\phi(\Delta)$) is determined from the shortest path length between vertices.  So (\ref{cosum}) will hold if the collection of paths going into $X_j$ is the same as the collection going out of $X_{j^*}.$  When such a relationship holds, we call $X_{j^*}$ the co-state of $X_j.$  Fig \ref{fcostate} illustrates such a pairing for our example problem, where $X_j=(2,1,4),$ with co-state $(4,1,2).$

Having now restated the issue as a graph question, we introduce the following definitions from graph theory:\cite{Chartrand}
\begin{itemize}
\item A digraph $D_1$ is {\it isomorphic} to digraph $D_2$ if there is a one-to-one mapping $\Phi$ from the vertex set of $D_1$ to the vertex set of $D_2$ such that $(u,v)$ is an edge of $D_1$ if and only if $(\Phi u,\Phi v)$ is an edge of $D_2.$
\item The {\it converse} of graph $D,$ (written $\overleftrightarrow{D}$) is the digraph with the same vertex set as $D,$ such that $(u,v)$ is an edge of $D$ if and only if $(v,u)$ is an edge of $\overleftrightarrow{D}$ (The converse is formed by reversing the direction of each edge.)
\item A digraph $D$ is {\it self-converse} if $\overleftrightarrow{D}$ is isomorphic to $D.$
\end{itemize}

\begin{lemma}\label{lemselfconverse} Let $G$ the digraph formed by taking configurations (unblocked states) as vertices and edges defined by the rearrangement resulting from a single Bernoulli failure $\delta_i.$  Then $G$ is self-converse. \end{lemma}
To prove the lemma, we demonstrate the self-converse structure by providing the necessary isomorphism.

\noindent {\bf Proof of lemma.} For vertex $X=(x_1, \dots, x_k)$ define map $\Phi$ (and notation $(\cdot)^*$by
\be
\label{ecostatedefn}
X^*:=\Phi(X) \equiv (x_k, \ldots, x_1).
\ee
$\Phi$ is clearly a one-to-one mapping of the vertices of $G$ onto itself.  Assume $(X,Y)$ is an edge of $G,$ (and, by definition, $(Y,X)$ is an edge in $\overleftrightarrow{G}.$)   Then $X+\delta_i=Y$ for some $i.$  Since vector addition in $\mathbb{R}^k$ is defined component-wise, reversing the order of components (by applying $\Phi$) does not alter the quality, so we have 
\be
\label{elem1}
X^*+\delta_i^*=Y^*.
\ee
By referring to the definition of $\delta_i,$ one easily sees that $\delta_i^*=-\delta_{k-i}.$ (Example: If $k=6,$ then $\delta_2=(0,1,-1,0,0,0),$ and $\delta_2^*=(0,0,0,-1,1,0)=-\delta_4.$) Substituting into \eqref{elem1} and rearranging, we have 
\[
Y^*+\delta_{k-i}=X^*,
\]
which implies there is an edge in $G$ from $Y^*$ to $X^*.$  So $(Y,X)$ an edge of $\overleftrightarrow{G}$ implies $(\Phi Y,\Phi X)$ is an edge of $G.$  Since each step above is reversible, $\Phi$ is the necessary isomorphism between $G$ and $\overleftrightarrow{G} . \qquad \Box$

\begin{figure}[htbp]
\includegraphics{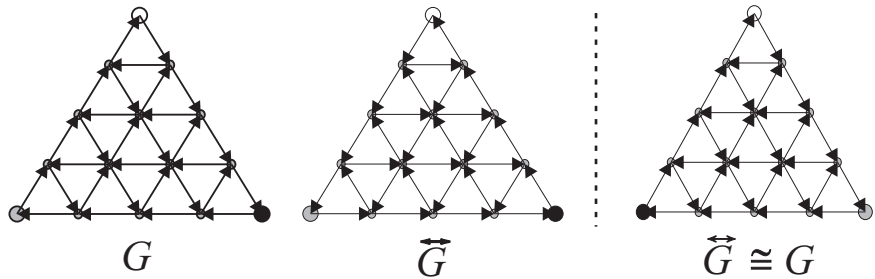} 
\caption{Using $k=3$ and $n=7,$ the first graph is $G,$ and the second is its converse, formed by reversing the direction of all the edges. The third graph is a simply a redrawing of the converse by reflecting the image over the hashed line.  The reflection preserves the direction of adjacency between vertices, so the graph is preserved. The new redrawing more clearly shows that $G$ is self-converse.}
\label{fconv}
\end{figure}

\vspace{.2in}
Fig \ref{fconv} illustrates the self-converse property of our example case.  

We remark that the $k$-fold symmetry of the graphs leads to $k$ different isomorphisms that prove the self-converse property.  However, we will consider (\ref{ecostatedefn}) as defining co-states.  The graphical structure {\it from} $X^*$ exactly matches the structure {\it into} $X,$ which ensures (\ref{cosum}) is satisfied and establishes the balance among the unblocked states.  Moreover, the self-converse structure tells us that for any two unblocked states $X$ and $Y,$
\[
p_{XY}=p_{Y^*X^*}
\]
However to prove the unblocked states are stationary, we must demonstrate an additional balance between the density transitioning into $X_j$ {\it from blocked states} with the density transitioning from $X_j^*$ {\it to blocked states.}

Let $W^b$ be a state with $b$ blockages and $W$ the unblocked state with the same physical configuration as $W^b.$  We recall that transition probabilities {\it from} $W^b$ are the same as the transition probabilities from $W$ (because the system does not have ``memory" of how it got to some configuration).  If the current density is $\nu,$ then the density of $W^b$ is $r^b/ \|\omega\|$ (by definition of $\nu$).  The density that transitions from blocked state $W^b$ to any unblocked state $X$ is 
\be
\label{elhs}
p_{WX}\frac{r^b}{\| \omega\|}.
\ee
This is exactly the same as the density which transitions from $X^*$ to $\Phi(W^b),$ computed as follows: the density of $X^*$ is $1/ \| \omega \|,$ so the density transition from $X^*$ to $W^*$ is $\frac{1}{\| \omega \|} p_{X^*W^*}.$  We apply (\ref{probr}) $b$ times to find that the density transitioned is
\be
\label{erhs}
\frac{1}{\| \omega \|} r^b p_{X^*W^*}=\frac{r^b}{\|\omega\|}p_{WX}.
\ee

This additional balance of transitioning density allows us to conclude that under assumed density $\nu,$ the unblocked state density is stationary.  To complete the proof, we invoke Lemma \ref{lemtoblocks} to claim that the blocked states are stationary as well.  Additionally, since $p_{ij}>0$ for all transition probabilities, it is known that the stationary density is unique\cite{Grimmett}.

\section{Calculating the Expected Value for Blockages. \label{sconclusion}}

Finding the stationary density for a Markov process enables us to answer a myriad of questions regarding the asymptotic behavior of the system.  In particular, we are now ready to tackle the question of the original problem statement.  Specifically, we seek to develop a formula for $B(k,n,r)=$ {\it the expected value of the fraction of the time steps for which a given worker is blocked.}  Since the problem is symmetric in each worker, without loss of generality, we analyze worker 1.  $\nu_i$ gives the expected value of the fraction of the time spent in state $i,$ so
\be
B(k,n,r)=\sum_{i \in \mathcal{A}} \nu_i,
\ee
where $\mathcal{A}$ is the index of all the states when worker 1 is blocked.  We note that
\be
\nu_i=\frac{\omega_i}{\|\omega\|} = \frac{r^{b_i}}{\|\omega\|},
\ee
where $b_i$ is the number of blockages in state $X_i.$  We first calculate $\|\omega\|$ by applying (\ref{nb}):
\be
\|\omega\| = \sum _ {b=0} ^ {k-1} N(b,k,n) r^b = \sum _ {b=0} ^ {k-1} {k \choose b}{n-b-1 \choose k-b-1} r^b.
\ee
We now consider the set $\mathcal{A}:$ When worker 1 is blocked, the state is given by 
\be
X=(1^*,x_2,\ldots,x_k).
\ee
So we may formally define
\be
\mathcal{A}=\left\{j \left| X_j=(1^*,x_2,\ldots,k_k) , \sum_{i=2}^k |x_i|=n-1 \right. \right\}.
\ee
We note that for any state $X_i, i \in \mathcal{A},$ the total number of blockages is exactly one more than the number of blockages in the 2nd through $k$th components.  We can use our previously determined counting formulas, but with these updated arguments, computing $N(b,k-1,n-1),$ where $b$ may range from $0$ to $k-2.$ We then find that the frequency of blockage of worker 1 is given by
\be
\label{etheformula}
B(k,n,r)  = \sum\limits_{i \in \mathcal{A}} \nu_i =  \sum\limits_{i \in \mathcal{A}} \frac{r^{b_i}}{\|\omega\|} = \frac{ \sum \limits_ {b=0} ^ {k-2} {k-1 \choose b}{n-b-2 \choose k-b-2} r^{b+1}}{\sum\limits _ {b=0} ^ {k-1} {k \choose b}{n-b-1 \choose k-b-1} r^b}. 
\ee
[Note: the exponent $r^{b+1}$ accounts for the formulation of one $1^*$ in the first component and $b$ $1^*s$ among components $x_2,\ldots,x_k.$]  

\section{Conclusion}

In this paper, we considered the problem of single-file motion of multiple agents on a closed path.  By using a moving reference frame of distance between workers (instead of actual bin location), we were able to simplify the geometry of the constraint.  In the $k-1$ dimensional state-space, the system reduces to a random walk with the coordinate planes as boundaries.  Distinct from the usual state transition graph, we proposed a decomposition digraph for transitions that allowed calculation of transition probabilities.  The self-converse structure of that digraph allowed us to identify the stationary density of the process.  We expect that the geometric approach may simplify similar boundary descriptions of other single file diffusion problems.  Additionally, we suspect that similar symmetries may be exploitable by decomposition of these processes into fundamental motion digraphs.

The original motivation for this problem was an actual operations analysis question involving zoning arrangements  to achieve maximum worker efficiency in warehousing systems.  Our simplified expression of that problem (as stated in this paper) implies that workers may move arbitrarily far around the warehouse in one time step before achieving a Bernoulli success.  We call this situation the ``infinite travel speed'' problem.  As a more realistic model, one would consider an additional constraint that a worker could have at most $f$ failures during a time step, since there must be some maximum speed of travel in a real system.  The closed form description of the resultant stationary density remains an open problem.

{}

\end{document}